\documentclass[11pt]{article}
\newcommand{\Mcb}{\mathcal{M}_\mathrm{cb}}
\renewcommand{\Theta}{\theta} 
\usepackage{amsmath,amssymb,amsfonts,diagrams}
\newenvironment{keywords}{\noindent\small {\it Keywords\/}:}{\vskip 4pt}
\newenvironment{classification}{\noindent\small 2000 {\it Mathematics Subject
Classification\/}:}{\vskip 12pt}

\renewcommand{\iff}{\quad\Longleftrightarrow\quad}

\newcommand{\comps}{{\mathbb C}}

\newcommand{\ints}{{\mathbb Z}}
\newcommand{\posints}{{\mathbb N}}

\newcommand{\torus}{{\mathbb T}}
\newcommand{\free}{{\mathbb F}}

\newcommand{\void}{\varnothing}
\newcommand{\tensor}{\otimes}
\newcommand{\ttensor}{\tilde{\otimes}}
\newcommand{\Tensor}{\hat{\otimes}}

\newcommand{\sdprod}{\rtimes}
\newcommand{\cstar}{{C^\ast}}

\newcommand{\re}{{\operatorname{Re}}}

\newcommand{\id}{{\mathrm{id}}}
\newcommand{\cb}{{\mathrm{cb}}}

\newcommand{\hull}{{\operatorname{hull}}}

\newcommand{\Hilbert}{{\mathfrak H}}

\newcommand{\M}{{\mathfrak M}}

\newcommand{\CB}{{\cal CB}}

\newcommand{\VN}{\operatorname{VN}}

\newcommand{\SL}{\operatorname{SL}}
\newcommand{\SIN}{\operatorname{SIN}}

\newcommand{\supp}{{\operatorname{supp}}}

\newcommand{\varcl}[1]{\overline{#1}}

\usepackage{amsthm,enumerate}
\theoremstyle{plain}
\newtheorem{theorem}{Theorem}[section]
\newtheorem{lemma}[theorem]{Lemma}
\newtheorem{corollary}[theorem]{Corollary}
\newtheorem{proposition}[theorem]{Proposition}
\theoremstyle{definition}
\newtheorem{definition}[theorem]{Definition}
\theoremstyle{remark}
\newtheorem*{remark}{Remark}
\newtheorem*{example}{Example}
\newtheorem*{rems}{Remarks}
\newtheorem*{exs}{Examples}
\newenvironment{remarks}{\begin{rems}\begin{enumerate}}{\end{enumerate}\end{rems}}

\newenvironment{items}{\begin{enumerate}[\rm (i)]}{\end{enumerate}}
\newenvironment{alphitems}{\begin{enumerate}[\rm (a)]}{\end{enumerate}}

\renewcommand{\baselinestretch}{1.2}
\title{Harmonic operators: the dual perspective}
\author{\textit{Matthias Neufang}\thanks{Research supported by NSERC and
the Mathematisches Forschungsinstitut Oberwolfach.} \and 
\textit{Volker Runde}\thanks{Research supported by NSERC and
the Mathematisches Forschungsinstitut Oberwolfach.}}
\date{}
\begin{document}
\maketitle
\begin{abstract}
The study of harmonic functions on a locally compact group $G$ has recently 
been transferred to a  ``non-commutative'' setting in two different directions:
C.-H.\ Chu and A.\ T.-M.\ Lau replaced the algebra $L^\infty(G)$ by the 
group von Neumann algebra $\VN(G)$ and the convolution action of a probability 
measure $\mu$ on $L^\infty(G)$ by the canonical action of a positive definite 
function $\sigma$ on $\VN(G)$; on the other hand, W.\ Jaworski and the 
first-named author replaced $L^\infty(G)$ by ${\cal B} (L^2(G))$ to which 
the convolution action by $\mu$ can be extended in a natural way. 
We establish a link between both approaches. The action of $\sigma$ 
on $\VN(G)$ can be extended to ${\cal B} (L^2(G))$. We study the 
corresponding space $\tilde{\cal H}_\sigma$ of ``$\sigma$-harmonic operators'',
i.e., fixed points in ${\cal B} (L^2(G))$ under the action of $\sigma$. 
We show, under mild conditions on either $\sigma$ or $G$, that 
$\tilde{\cal H}_\sigma$ is in fact a von Neumann subalgebra of 
${\cal B} (L^2(G))$. Our investigation of $\tilde{\cal H}_\sigma$ relies, in 
particular, on a notion of support for an arbitrary operator 
in ${\cal B} (L^2(G))$ that extends Eymard's definition for elements of 
$\VN(G)$. Finally, we present an approach to $\tilde{\cal H}_\sigma$ via 
ideals in ${\cal T} (L^2(G))$---where ${\cal T}(L^2(G))$ denotes the trace
class operators on $L^2(G)$, but equipped with a product different from
composition---, as it was pioneered for harmonic functions by G.\ A.\ Willis.
\end{abstract}
\begin{sloppy}
\end{sloppy}
\section*{Introduction}
Let $G$ be a locally compact group, and let $\mu$ be a probability measure on
$G$ whose support generates $G$. A function $\phi \in L^\infty(G)$ is called
\emph{$\mu$-harmonic} if $\mu \ast \phi = \phi$. Harmonic functions play
a crucial r\^ole for the investigation of random walks on locally compact 
groups and are extensively studied (see \cite{Aze} and \cite{Fur}, for
example). The collection of all $\mu$-harmonic functions is denoted by
${\cal H}_\mu$; it is easily seen to be a $w^\ast$-closed subspace of
$L^\infty(G)$, but is a subalgebra only if it consists of the constants alone.
Nevertheless, ${\cal H}_\mu$ can be equipped with a product---different,
of course, from the one in $L^\infty(G)$---turning it into an abelian von
Neumann algebra (\cite{Aze}).
\par
Recently, the notion of harmonicity has been ``quantized'' in two, seemingly
entirely different directions.
\par
One such quantization was introduced and studied by C.-H.\
Chu and A.\ T.-M.\ Lau in \cite{CL}. Their approach can be considered dual to 
the classical concept of a harmonic function. In \cite{Eym}, P.\ Eymard
introduced the so-called Fourier algebra $A(G)$ and Fourier--Stieltjes algebra
$B(G)$ of $G$. If $G$ is abelian with dual group $\hat{G}$, the Fourier and
Fourier--Stieltjes transforms, respectively, yield isometric isomorphisms
$A(G) \cong L^1(\hat{G})$ and $B(G) \cong M(\hat{G})$. Instead of looking
at harmonicity with respect to a probability measure, Chu and Lau consider
harmonicity of functionals on $A(G)$, i.e., of elements of the group von
Neumann algebra $\VN(G)$, with respect to a positive definite function 
$\sigma \in B(G)$: an operator $T \in \VN(G)$ is said to be $\sigma$-harmonic 
if $\sigma \cdot T = T$, where $\cdot$ is the canonical module action of
$B(G)$ on $\VN(G)$. The collection of all $\sigma$-harmonic functionals is
denoted by ${\cal H}_\sigma$. Even though, this new theory seems to parallel 
the classical theory of harmonic functions on the surface, it is, in fact,
strikingly different. For instance, ${\cal H}_\sigma$ is \emph{always} a
von Neumann \emph{subalgebra} of $\VN(G)$ (\cite[Remark 3.2.11]{CL}).
\par
A completely different type of quantization was recently carried out by
W.\ Jaworski and the first-named author (\cite{JN}). Starting point is the
the result by F.\ Ghahramani (\cite{Gha}) that there is a natural isometric 
representation $\Theta$ of the measure algebra $M(G)$ on ${\cal B}(L^2(G))$,
such that for $\mu \in M(G)$ and $\phi \in L^\infty(G)$---viewed as a 
multiplication operator on $L^2(G)$---we have $\Theta(\mu)\phi = 
\mu \ast \phi$. Hence, the authors of \cite{JN} define an operator
$T \in {\cal B}(L^2(G))$ to be $\mu$-harmonic for a probability measure
$\mu$ if $\Theta(\mu)(T) = T$. The collection of all $\mu$-harmonic operators
is denoted by $\tilde{\cal H}_\mu$. Like ${\cal H}_\mu$, the space
$\tilde{\cal H}_\mu$ is a von Neumann algebra, but with a product usually
different from the one in ${\cal B}(L^2(G))$;
in fact, $\tilde{\cal H}_\mu$ can be described as the crossed product
of ${\cal H}_\mu$ with $G$, where the action of $G$ on ${\cal H}_\mu$ is
given by left translation (\cite[Proposition 6.3]{JN}). 
\par
In the present paper, we extend Chu's and Lau's notion of $\sigma$-harmonicity from $\VN(G)$ to ${\cal B}(L^2(G))$ in a way that parallels the extension of $\mu$-harmonicity from $L^\infty(G)$ to ${\cal B}(L^2(G))$ in \cite{JN}.
\par
As the predual space of the operator algebra $\VN(G)$, the Fourier algebra
carries a canonical operator space structure. (For the theory of operator
spaces, we refer to \cite{ER}.) By $\Mcb(A(G))$, we denote the
completely bounded multipliers of $A(G)$, i.e., those functions on $G$
that induce completely bounded multiplication operators on $A(G)$.
Obviously, $\Mcb(A(G))$ is a commutative Banach algebra, and it
contains $B(G)$ (with equality if and only if $G$ is amenable).
In \cite{NRS}, the first-named author, Z.-J.\ Ruan, and N.\ Spronk
constructed a representation $\hat{\Theta}$ of $\Mcb(A(G))$ on
${\cal B}(L^2(G))$ which displays a perfect duality with Ghahramani's
representation of $M(G)$ (for details, see \cite{NRS}); in particular, it
extends the action of $B(G)$ on $\VN(G)$ to ${\cal B}(L^2(G))$. For $\sigma
\in \Mcb(A(G))$, it is then natural to define $T \in {\cal B}(L^2(G))$ to be
$\sigma$-harmonic if $\hat{\Theta}(\sigma)(T) = T$. We denote the collection
of all $\sigma$-harmonic operators by $\tilde{\cal H}_\sigma$.
\par
The aim of this paper is to explore the connections between this setting
and the two quantizations from \cite{CL} and \cite{JN}. For instance,
one of our main results is that---under very mild hypotheses which are 
always satisfied if $G$ is amenable or the free group in two 
generators---$\tilde{\cal H}_\sigma$ is a von 
Neumann subalgebra of ${\cal B}(L^2(G))$, namely the von Neumann subalgebra 
of ${\cal B}(L^2(G))$ generated by ${\cal H}_\sigma$ and $L^\infty(G)$.
\par
The paper is organized as follows. 
\par
First, we fix our notation and terminology  
and review a few basic facts about harmonic functions (see \cite{JN}) 
and harmonic functionals in $\VN(G)$ (as introduced and studied in \cite{CL}).
We also recall results from \cite{NRS} on the representation of $\Mcb(A(G))$ 
on ${\cal B}(L^2(G))$. 
\par
Section \ref{adapted} gives a characterization of 
adaptedness in terms of the Fourier--Stieltjes transform (in the framework of 
locally compact \emph{abelian} groups). This motivates our definition of the 
dual notion of adaptedness for positive definite functions. 
\par
In Section \ref{support}, we introduce the support of an arbitrary 
(bounded linear) operator on $L^2(G)$ such that it generalizes Eymard's 
corresponding notion for elements of $\VN(G)$. In order to prepare the ground
for our main results, we investigate the properties of our notion of support
in detail. 
\par
Those main results are contained in Section \ref{dualpic}. Using 
the representation of $\Mcb(A(G))$ from \cite{NRS}, we naturally extend the 
notion of a $\sigma$-harmonic functional 
in $\VN (G)$ (\cite{CL}) to the one of a $\sigma$-harmonic operator in 
${\cal B}(L^2(G))$, where $\sigma$ is a positive definite function on $G$.
Assuming either a very mild condition on $G$---the approximation property
(\cite{HK})---or that $\sigma$ belongs to $A(G)$, we show that the space of 
$\sigma$-harmonic
operators $\tilde{\cal H}_\sigma$ is always a von Neumann subalgebra of 
${\cal B}(L^2(G))$ and, in fact, precisely the von Neumann algebra 
generated by $L^\infty(G)$ and the algebra ${\cal H}_\sigma$ of 
$\sigma$-harmonic functionals in $\VN (G)$. This may be viewed as a result 
``dual'' to the corresponding characterization (see \cite{JN}) 
of the von Neumann algebra of $\mu$-harmonic operators in ${\cal B}(L^2(G))$ 
as a crossed product over the algebra of $\mu$-harmonic functions in 
$L^\infty(G)$, where $\mu$ is a probability measure on $G$. 
Our notion of support for an arbitrary operator on $L^2(G)$ allows for
an alternative characterization of $\tilde{\cal H}_\sigma$---at least if $G$
has the approximation property---, namely as the operators in 
${\cal B}(L^2(G))$ with support contained in the subgroup 
$G_\sigma = \sigma^{-1} (\{ 1 \})$ of $G$. 
\par
Finally, in Section \ref{ideals}, we present an approach to harmonic operators 
in ${\cal B}(L^2(G))$ via ideals in the predual ${\cal T}(L^2(G))$ in the 
spirit of \cite{Wil} (see also \cite[Section 3]{JN}). This makes it 
necessary to endow the space ${\cal T}(L^2(G))$ with a product very different 
from the composition of operators; this product arises naturally when 
when one regards ${\cal B}(L^2(G))$ as a Hopf--von Neumann algebra, with a 
co-multiplication naturally extending the one of 
$\VN(G)$ (see \cite{NRS}, \cite{PR}, and \cite{JN}). 
\subsubsection*{Acknowledgments}
Most of this paper was written during an RiP ($=$ ``Research in Paradise'')
stay by the authors at the Mathematisches Forschungsinstitut Oberwolfach in
May 2005. We acknowledge this support and the uniquely stimulating atmosphere
at Oberwolfach with gratitude.
\par
We would also like to thank Tony Lau for helping with the proof of Lemma
\ref{wdense}---in particular, for bringing \cite{Lau} to our attention---and
Zhiguo Hu for catching some oversights in an earlier version of this paper. Finally, thanks are due to the referee for his helpful suggestions.
\section{Preliminaries}
Throughout this section, let $G$ denote a locally compact group equipped with
left Haar measure. 
We follow standard notation and terminology of abstract harmonic analysis as, e.g., in \cite{HR}. 
For $p \in [1,\infty]$, we denote the corresponding $L^p$-space with respect to left Haar measure on $G$ by $L^p(G)$. By $M(G)$,
we denote the algebra of all complex, regular Borel measures; $M^1(G)$
stands for the probability measures in $M(G)$. We note that $L^\infty(G)$
is a von Neumann algebra acting naturally on $L^2(G)$ as multiplication
operators: for $\phi \in L^\infty(G)$, let $M_\phi \in {\cal B}(L^2(G))$
be given by $(M_\phi \xi)(x) = \phi(x) \xi(x)$ for $\xi \in L^2(G)$ and
$x \in G$. (For the sake of notational simplicity, we shall sometimes write
$\phi$ and $M_\phi$ interchangeably.) 
\par
The Banach algebra $M(G)$ acts on $L^p(G)$ for $p \in [1,\infty]$ via 
convolution from the left. For $\mu \in M(G)$, we define 
$\phi \in L^\infty(G)$ as \emph{$\mu$-harmonic} if $\mu \ast \phi = \phi$
and set
\[
  {\cal H}_\mu := \{ \phi \in L^\infty(G) : \text{$\phi$ is $\mu$-harmonic} \}.
\]
(For the motivation for the name ``$\mu$-harmonic function'', see the
introduction of \cite{CL}). Usually, $\mu$-harmonic functions are considered 
only for $\mu \in M^1(G)$. To avoid pathologies, we shall assume 
$\mu$ to be \emph{adapted}, i.e., $\langle \supp \mu \rangle$, the closed 
subgroup of $G$ generated by the support of $\mu$, is all of $G$. For abelian $G$ and adapted $\mu \in M^1(G)$, the classical Choquet--Deny theorem asserts that ${\cal H}_\mu$ consists only of the constant functions. For a general locally compact group $G$, a measure $\mu \in M^1(G)$ such that ${\cal H}_\mu \cong \comps$ exists if and only if $G$ is amenable and $\sigma$-compact (\cite[Proposition 2.1.3]{CL}). Even though ${\cal H}_\mu$ is not a von Neumann subalgebra of $L^\infty(G)$---except if ${\cal H}_\mu \cong \comps$---, there is a natural product on ${\cal H}_\mu$ turning it into an abelian von Neumann algebra in its own right: given $f,g \in {\cal H}_\mu$, one can show that the limit 
\[
  \text{$w^\ast$-}\lim_n \int_G \rho(x) (fg) \, d\mu^{\ast n} (x)
\]
exists in $L^\infty(G)$ and is an explicit formula for the multiplication
in ${\cal H}_\mu$ (here, $\mu^{\ast n}$ stands for the $n$-th convolution 
power of $\mu$, and $\rho(x)$ denotes right 
translation by $x$, i.e., $(\rho(x)\phi)(y) :=\phi(yx)$ for 
$\phi \!: G \to \comps$ and $x,y \in G$). 
For the classical theory of random walks and 
harmonic functions, see \cite{Aze}, \cite{Nev}, and \cite{Rev}, for example.
\par
For any $\phi \!: G \to \comps$ and $x \in G$, we define $\lambda(x)\phi
\!: G \to \comps$ by letting $(\lambda(x)\phi)(y) := \phi(x^{-1}y)$ for
$x \in G$. For fixed $x \in G$, the map $L^2(G) \ni \xi \mapsto \lambda(x)\xi$
is a unitary operator on $L^2(G)$, and 
\[
  \lambda \!: G \to {\cal B}(L^2(G)), \quad x \mapsto \lambda(x)
\]
is a unitary representation of $G$ on $L^2(G)$, the \emph{left regular 
representation} of $G$. The \emph{group von Neumann algebra} $\VN(G)$ of
$G$ is defined as $\VN(G) := \lambda(G)''$. The \emph{Fourier algebra}
$A(G)$ of $G$---introduced by P.\ Eymard in \cite{Eym}---is the (unique) 
predual of $\VN(G)$. 
\par 
As any operator algebra, $\VN(G)$ is a concrete operator space. (For the 
theory of operator spaces, our reference is \cite{ER}). Since the category of 
operator spaces allows for a natural duality theory, the dual space 
$\VN(G)^\ast$---and thus its subspace $A(G)$---is an operator space in a 
canonical manner. In particular, we may speak of completely bounded maps
on $A(G)$. (Following \cite{ER}, we denote the completely bounded maps
on an operator space $E$ by $\CB(E)$.) A \emph{multiplier} of $A(G)$ is a
function $\sigma$ on $G$ such that $\sigma A(G) \subset A(G)$; a multiplier is 
completely bounded if the corresponding multiplication operator is 
completely bounded. The collection of all completely bounded multipliers
is denoted by $\Mcb(A(G))$. It naturally inherits an operator space structure
from $\CB(A(G))$ which turns it into a commutative, completely contractive
Banach algebra, i.e., a Banach algebra which is an operator space such that
multiplication is completely contractive. An important subalgebra of
$\Mcb(A(G))$ is the \emph{Fourier--Stieltjes algebra} $B(G)$ of $G$ (introduced also in \cite{Eym}). 
It consists of all functions of the form $G \ni x \mapsto \langle \pi(x) \xi, \eta \rangle$, 
where $\pi$ is a (always strongly continuous) unitary representation of $G$ on a Hilbert space $\Hilbert$ and
$\xi, \eta \in \Hilbert$. It can be identified with the dual space of 
the full group $\cstar$-algebra $\cstar(G)$, i.e., of the enveloping 
$\cstar$-algebra of $L^1(G)$, and thus also has a canonical operator space
structure. It contains $A(G)$ as a closed ideal and thus canonically embeds 
into $\Mcb(A(G))$; this embedding is always completely contractive, but is
an isometric isomorphism if and only if $G$ is amenable. For more details
and references to the original literature, see \cite{Spr}.
\par
It is immediate from the definition of $\Mcb(A(G))$ 
that $A(G)$ is a completely contractive Banach $\Mcb(A(G))$-module through
pointwise multiplication. Consequently, $\VN(G)$ carries a dual 
$\Mcb(A(G))$-module structure, namely
\[
  \langle \phi, \sigma \cdot T \rangle := \langle \phi \sigma, T \rangle
  \qquad (\sigma \in \Mcb(A(G)), \, T \in \VN(G), \, \phi \in A(G)).
\]
Slightly generalizing the definition from \cite{CL}, we say that, for
$\sigma \in \Mcb(A(G))$, a von Neumann algebra element $T \in V(G)$, is
\emph{$\sigma$-harmonic} if $\sigma \cdot T = T$, and define
\[
  {\cal H}_\sigma := \{ T \in \VN(G) : \text{$T$ is $\sigma$-harmonic} \}.
\]
(Chu and Lau consider ${\cal H}_\sigma$ only for 
$\sigma \in P^1(G)$, which denotes the positive definite functions of norm one
in $B(G)$.) For abelian $G$ with dual group $\hat{G}$, the Fourier and 
Fourier--Stieltjes transforms, respectively, yield isometric isomorphism
$A(G) \cong L^1(\hat{G})$ and $B(G) \cong M(\hat{G})$, and conjugation with
the Plancherel transform yields that
\begin{equation} \label{planschi}
  {\cal H}_\mu \cong {\cal H}_{\hat{\mu}} \qquad (\mu \in M(G)).
\end{equation}
Despite the formal analogies with harmonic functions, the harmonic
functionals according to Chu and Lau display a strikingly different behavior:
for $\sigma \in P^1(G)$, the set
\[
  G_\sigma := \{ x \in G : \sigma(x) = 1 \}
\]
is a closed subgroup of $G$ such that ${\cal H}_\sigma = \lambda(G_\sigma)''$
(\cite[Proposition 3.2.10]{CL}); in particular, ${\cal H}_\sigma$ is
a von Neumann subalgebra of ${\cal B}(L^2(G))$.
\par
Finally, we require a construction from \cite{NRS}: there the first-named 
author, Z.-J.\ Ruan, and N.\ Spronk
define complete isometry $\hat{\Theta} \!: \Mcb(A(G)) \to 
\CB({\cal B}(L^2(G)))$ with the following properties:
\begin{itemize}
\item $\hat{\Theta}$ is a unital, $w^\ast$-$w^\ast$-continuous algebra
homomorphism;
\item $\hat{\Theta}(\Mcb(A(G)))$ consists precisely of those normal operators
in $\CB({\cal B}(L^2(G)))$ which are $L^\infty(G)$-bimodule homomorphisms and leave
$\VN(G)$ invariant;
\item for $\phi \in \Mcb(A(G))$ and $T \in \VN(G)$, we have
$\hat{\Theta}(\phi)(T) = \phi \cdot T$, i.e., the action of $\Mcb(A(G))$
on ${\cal B}(L^2(G))$ induced by $\hat{\Theta}$ extends the canonical one
of $\Mcb(A(G))$ on $\VN(G)$.
\end{itemize}
Under many aspects, $\hat{\Theta}$ can be viewed as dual to Ghahramani's
representation of $M(G)$ on ${\cal B} (L^2(G))$ (for details, see Section \ref{dualpic}, and \cite{NRS}).
\par
For later use, we indicate how $\hat{\Theta}$ is defined. Given $\phi
\in \Mcb(A(G))$, there are a Hilbert space $\Hilbert$ and continuous
functions $\boldsymbol{\xi}, \boldsymbol{\eta} \!: G \to \Hilbert$ such that
\[
  \phi(x^{-1}y) = \langle \boldsymbol{\xi}(y), 
                  \boldsymbol{\eta}(x) \rangle \qquad (x,y \in G)
\]
(see \cite{Jol} for an accessible proof). Let $( e_i )_{i \in \mathbb I}$
be an orthonormal basis, and define
\[
  \phi_i(x) := \langle e_i, \boldsymbol{\eta}(x) \rangle
  \quad\text{and}\quad
  \psi_i(x) := \langle \boldsymbol{\xi}(x), e_i \rangle
  \qquad (i \in \mathbb{I}, \, x \in G).
\]
Then we have
\[
  \hat{\Theta}(\phi)(T) := \sum_{i \in \mathbb I} M_{\phi_i(x)} T M_{\psi_i(x)}.
\]
For more details, see \cite{NRS}.
\section{A dual characterization of adapted probability measures} \label{adapted}
Let $G$ be a locally compact group, and let $\mu \in M^1(G)$. It is fair to
say that almost nothing can be said about ${\cal H}_\mu$ unless
$\mu$ is adapted. In order to have a dual notion of adaptedness, i.e.,
in the context of harmonic functionals as in \cite{CL}, we first
characterize the adapted probability measures on a locally
compact \emph{abelian} group in terms of their Fourier--Stieltjes transforms.
\begin{proposition} \label{adapprop}
Let $G$ be a locally compact abelian group. Then the following are equivalent
for $\mu \in M^1(G)$:
\begin{items}
\item $\mu$ is adapted;
\item if $\gamma \in \hat{G}$ is such that $\gamma |_{\supp \, \mu} \equiv 1$,
then $\gamma = 1$;
\item $\{ \gamma \in \hat{G} : \hat{\mu}(\gamma) = 1 \} = \{ 1 \}$.
\end{items}
\end{proposition}
\begin{proof}
(i) $\Longleftrightarrow$ (ii) and (iii) $\Longrightarrow$ (i) are straightforward.
\par 
(ii) $\Longrightarrow$ (iii): Let $\gamma \in \hat{G}$ be such that
\[
  \hat{\mu}(\gamma) = \int_G \overline{\gamma(x)} \, d\mu(x) = 1,
\]
so that
\[
  1 = \overline{\hat{\mu}(\gamma)} = \overline{\int_G 
  \overline{\gamma(x)} \, d\mu(x)} =
  \int_G \gamma(x) \, d\mu(x) 
\]
and thus
\[
  \int_G \re \, \gamma(x) \,d \mu(x) = \frac{1}{2} \left(
  \int_G \overline{\gamma(x)} \, d\mu(x) + \int_G \gamma(x) \, d\mu(x)
  \right) = 1.
\]
Since $\gamma(G) \subset \torus$, we have $(\re \, \gamma)(G) \subset [-1,1]$,
so that
\[
  0 = \int_G 1 \, d\mu(x) - \int_G \re \, \gamma(x) \, d\mu(x) =
  \int_{\supp \, \mu} \underbrace{(1 - \re \, \gamma(x))}_{\geq 0} \, d\mu(x).
\]
By the continuity of $\gamma$, this means that $\re \, \gamma |_{\supp \, \mu} \equiv 1$, 
thus $\gamma |_{\supp \, \mu} \equiv 1$, and therefore $\gamma = 1$ by (ii). 
\end{proof}
\par
Using (\ref{planschi}) and \cite[Proposition 3.2.10]{CL}, we obtain a 
dual approach to the Choquet--Deny theorem:
\begin{corollary}
Let $G$ be a locally compact abelian group, and let $\mu \in M^1(G)$ be
adapted. Then ${\cal H}_\mu \cong \comps$ holds.
\end{corollary}
\par
In view of Proposition \ref{adapprop}, we 
define for general locally compact groups:
\begin{definition} \label{adapdef}
Let $G$ be a locally compact group. Then we call $\sigma \in P^1(G)$ 
\emph{adapted} if $G_\sigma = \{ e \}$.
\end{definition}
\begin{remark}
In terms of Definition \ref{adapdef}, Proposition \ref{adapprop} can be
reformulated as follows: a probability measure on a locally compact abelian group
is adapted if and only if its Fourier--Stieltjes transform is adapted.
\end{remark}
\par
It is well known that there is an adapted probability measure on a locally
compact group $G$ if and only if $G$ is $\sigma$-compact. Since a locally
compact abelian group is $\sigma$-compact if and only if its dual is
first countable, it is immediate that there is an adapted positive definite
function on such a group if and only if it is first countable.
\par
The following proposition extends this to general locally compact groups:
\begin{proposition} \label{adapprop2}
The following are equivalent for a locally compact group:
\begin{items}
\item there is an adapted $\sigma \in P^1(G)$;
\item $G$ is first countable;
\item there is an adapted $\sigma \in A(G) \cap P^1(G)$.
\end{items}
\end{proposition}
\begin{proof}
(i) $\Longrightarrow$ (ii): Let $\sigma \in P^1(G)$ be adapted, fix a compact
neighborhood $U$ of $e$, and define,
for $n \in \posints$, 
\[
  U_n := \left\{ x \in U : | \sigma(x) - 1 | < \frac{1}{n} \right\}.
\] 
From the continuity of $\sigma$, it is clear that $\{ U_n : n \in \posints \}$
consists of neighborhoods of $e$. Let $V$ be a neighborhood
of $e$, and suppose without loss of generality that $V$ is open and contained
in $U$. Since $\sigma$ is continuous and adapted, and since $U \setminus V$
is compact, we have that
\[
  \epsilon_0 := \inf \{ | \sigma(x) - 1 | : x \in U \setminus V \} > 0.
\]
Choose $n_0 \in \posints$ so large that $\frac{1}{n_0} \leq \epsilon_0$.
It follows that $U_{n_0} \subset V$. Consequently, $\{ U_n : n \in \posints \}$
is a base of neighborhoods of $e$.
\par
Through translation, we see that every point of $G$ has a countable 
base of neighborhoods. 
\par
(ii) $\Longrightarrow$ (iii): Let $\{ U_n : n \in \posints \}$ be a base
of neighborhoods of $e$, and suppose without loss of generality that
$U_{n+1} \subset U_n$ for $n \in \posints$. For each $n \in \posints$, 
there is $\sigma_n \in A(G) \cap P^1(G)$ with $\supp \, \sigma_n \subset U_n$.
Define
\[
  \sigma := \sum_{n=1}^\infty \frac{1}{2^n} \sigma_n
\]
so that, clearly, $\sigma \in A(G) \cap P^1(G)$. Let $x \in G \setminus \{ e \}$.
Since $\{ U_n : n \in \posints \}$ is a base of neighborhoods of $e$, there
is $n_0 \in \posints$ such that $x \notin U_{n_0}$ and thus $x \notin U_n$ 
for $n \geq n_0$. It follows that
\[
  |\sigma(x)| \leq \sum_{n=1}^{n_0-1} \frac{1}{2^n} |\sigma_n(x)| 
  \leq  \sum_{n=1}^{n_0-1} \frac{1}{2^n} < 1.
\]
This proves (iii).
\par
Finally, (iii) $\Longrightarrow$ (i) is trivial.
\end{proof}
\section{The support of an operator on $L^2(G)$} \label{support} 
In \cite{Eym}, P.\ Eymard introduced the notion of support for elements 
of a group von Neumann algebra: if $G$ is a locally compact group and $T
\in \VN(G)$, then the \emph{support} $\supp \, T$ of $T$ is defined to
consist of those $x \in G$ such that $\phi(x) = 0$ for all $\phi \in A(G)$
with $\phi \cdot T = 0$. 
\par
Using $\hat{\Theta}$, this notion can naturally be extended to arbitrary
operators on $L^2(G)$:
\begin{definition} \label{suppdef}
Let $G$ be a locally compact group, and let $T \in {\cal B}(L^2(G))$.
Then the \emph{support} $\supp \, T$ of $T$ is defined to be
\[
  \supp \, T := \{ x \in G: \text{$\phi(x) = 0$ for all $\phi \in A(G)$
  with $\hat{\Theta}(\phi)(T) = 0$} \}.
\]
\end{definition}
\begin{remarks}
\item For operators in $\VN(G)$, this notion of support coincides with the
one from \cite{Eym}.
\item The support of an operator is obviously a closed subset of $G$. 
\end{remarks}
\par
We first prove a few general assertions on the support of an operator.
\begin{proposition} \label{supportgeneral}
Let $G$ be a locally compact group. Then, for $\phi \in A(G)$ and $T \in {\cal B}(L^2(G))$, we have 
\[
  \supp \, \hat{\Theta}(\phi)(T) \subseteq \supp \, \phi \cap \supp \, T.
\] 
\end{proposition}
\begin{proof}
Fix $\phi \in A(G)$ and $T \in {\cal B}(L^2(G))$. 
\par
Since $\hat{\Theta}$ is multiplicative and $A(G)$ is commutative, it is 
straightforward that any $\psi \in A(G)$ with $\hat{\Theta}(\psi)(T) = 0$ 
satisfies
$\hat{\Theta}(\psi)(\hat{\Theta}(\phi)(T)) = 0$ as well. From Definition 
\ref{supportgeneral}, it is then immediate that $\supp \, \hat{\Theta}(\phi)(T)
\subset \supp \, T$.
\par
To see that $\supp \, \hat{\Theta}(\phi)(T) \subset \supp \, \phi$ as well, let
$x \in G \setminus \supp \, \phi$, and assume towards a 
contradiction that $x \in \supp \, \hat{\Theta}(\phi)(T)$. 
Since $x \notin \supp \, \phi$, there is a neighborhood $V$ of $x$ such that 
$\phi |_V \equiv 0$. Take $\psi \in A(G)$ such that $\supp \, \psi \subseteq V$
and $\psi (x) = 1$. Then $\psi \phi \equiv 0$ holds, so that 
$\hat{\Theta}(\psi) (\hat{\Theta}(\phi)(T)) = 
\hat{\Theta}(\psi \phi) (T) = 0$. 
Since $x \in \supp \, \hat{\Theta}(\phi)(T)$, this implies that 
$\psi (x) = 0$, which is a contradiction. 
\end{proof} 
\par
The following is an alternative description of the support
of an operator:
\begin{proposition} \label{suppideal}
Let $G$ be a locally compact group, and let $T \in {\cal B}(L^2(G))$.
Then 
\begin{equation} \label{ideal}
  \{ \phi \in A(G) : \hat{\Theta}(\phi)(T) = 0 \}
\end{equation}
is an ideal of $A(G)$ whose hull equals $\supp \, T$.
\end{proposition}
\begin{proof}
Let the ideal (\ref{ideal}) be denoted by $I$. Then we have for $x \in G$ that
\begin{align*}
  x \notin \supp \, T & \iff \text{there is $\phi \in A(G)$ such that $\phi(x) \neq 0$ and
  $\hat{\Theta}(\phi)(T) = 0$} \\
  & \iff \text{there is $\phi \in I$ with $\phi(x) \neq 0$} \\
  & \iff x \notin \hull(I).
\end{align*}
This completes the proof.
\end{proof}
\par
Since $A(G)$ is Tauberian for any locally compact group $G$, the following
is clear:
\begin{corollary} \label{suppcor0}
Let $G$ be a locally compact group, and let $T \in {\cal B}(L^2(G))$ be
such that $\supp \, T = \void$. Then $\hat{\Theta}(\phi)(T) = 0$ holds
for all $\phi \in A(G)$.
\end{corollary}
\par
Recall (\cite{HK}) that a locally compact group $G$ is said to have the
\emph{approximation property} if the constant function $1$ lies in the
$w^\ast$-closure of $A(G)$ in $\Mcb(A(G))$. Clearly, every amenable, locally
compact group has the approximation property (by Leptin's theorem), but so
does also every weakly amenable group in the sense of \cite{dCH}, such as
$\free_2$, the free group in two generators. Nevertheless, the approximation
property is weaker than weak amenability: by \cite[Corollary 1.17 and
Remark 3.10]{HK}, the group $\ints^2 \sdprod \SL(2,\ints)$ has the
approximation property, but is not weakly amenable.
\begin{proposition} \label{supprop0}
Let $G$ be a locally compact group with the approximation property. Then
the following are equivalent for $T \in {\cal B}(L^2(G))$:
\begin{items}
\item $\supp \, T = \void$;
\item $T = 0$.
\end{items}
\end{proposition}
\begin{proof}
Of course, only (i) $\Longrightarrow$ (ii) needs proof.
\par
Let $(e_\alpha)_\alpha$ be a net in $A(G)$ that converges to $1$ in the
$w^\ast$-topology of $\Mcb(A(G))$. Since $\hat{\Theta}$ is unital and
$w^\ast$-$w^\ast$-continuous (\cite[Theorem 4.5]{NRS}), it follows that
$\id_{{\cal B}(L^2(G))} = \text{$w^\ast$-}\lim_\alpha \hat{\Theta}(e_\alpha)$
and thus, by Corollary \ref{suppcor0},
\[
  T = \hat{\Theta}(1)(T) = 
  \text{$w^\ast$-}\lim_\alpha \hat{\Theta}(e_\alpha)(T) = 0.
\]
This proves (ii).
\end{proof}
\begin{remark}
It is well possible that Proposition \ref{supprop0} is true for \emph{every}
locally compact group.
\end{remark}
\par
Let $G$ be a locally compact group, and let $F \subset G$ be closed. As
is customary, we write
\[
  I(F) := \{ \phi \in A(G) : \phi |_F \equiv 0 \}.
\]
We also define
\[
  {\cal B}_F(L^2(G)) := \{ T \in {\cal B}(L^2(G)) : \supp \, T \subset F \}.
\]
\begin{lemma} \label{supplem2}
Let $G$ be a locally compact group, and let $F \subset G$ be a set of
synthesis for $A(G)$. Then $\hat{\Theta}(\phi)(T) = 0$ holds for all
$\phi \in I(F)$ and $T \in {\cal B}_F(L^2(G))$.
\end{lemma}
\begin{proof}
Let $T \in {\cal B}_F(L^2(G))$, and set
\[
  I := \{ \phi \in A(G) : \hat{\Theta}(\phi)(T) = 0 \}.
\] 
By Proposition \ref{suppideal}, we have
\[
  \hull(I) = \supp \, T \subset F.
\]
Since $F$ is a set of synthesis, this means that $I(F) \subset I$.
\end{proof}
\par
Since the multiplication of $\Mcb(A(G))$ is separately $w^\ast$-continuous
(see \cite{FRS}, for instance), we obtain immediately (from the 
$w^\ast$-$w^\ast$-continuity of $\hat{\Theta}$):
\begin{corollary} \label{suppcor1}
Let $G$ be a locally compact group with the approximation property,
let $F \subset G$ be a set of synthesis, and let $T \in {\cal B}_F(L^2(G))$. 
Then $\hat{\Theta}(\phi)(T) =0$ for all $\phi \in \Mcb(A(G))$ with 
$\phi |_F \equiv 0$.
\end{corollary}
\begin{proof}
Let $\phi \in \Mcb(A(G))$ be such that $\phi |_F \equiv 0$, and let 
$( e_\alpha )_\alpha$ be a net in $A(G)$ converging to $1$ in the 
$w^\ast$-topology of $\Mcb(A(G))$. Then $( e_\alpha \phi )_\alpha$ is 
a net in $I(F)$ that is $w^\ast$-convergent to $\phi$. Lemma \ref{supplem2}
and the $w^\ast$-$w^\ast$-continuity of $\hat{\Theta}$ then yield the claim.
\end{proof}
\par
Let $G$ be a locally compact group, let $F \subset G$ be a set of synthesis,
and let $T \in {\cal B}_F(L^2(G))$.
Then it is clear from Lemma \ref{supplem2} that $\hat{\Theta}(\phi)(T)$
for $\phi \in A(G)$ depends only on $\phi |_F$ (if $G$ has the approximation
property, this is even true for all $\phi \in \Mcb(A(G))$ by
Corollary \ref{suppcor1}).
\begin{proposition} \label{supprop1}
Let $G$ be a locally compact group, let $H$ be a closed subgroup of $G$,
and let $T \in {\cal B}_H(L^2(G))$. Then:
\begin{items}
\item if there is $\sigma \in A(G)$ with $\sigma |_H = 1$ and
$\hat{\Theta}(\sigma)(T) = T$, then $\hat{\Theta}(\phi)(T) = \phi(e) T$
for all $\phi \in A(G)$ which are constant on $H$;
\item if $G$ has the approximation property, then 
$\hat{\Theta}(\phi)(T) = \phi(e) T$ holds for all $\phi \in \Mcb(A(G))$
which are constant on $H$.
\end{items}
\end{proposition}
\begin{proof}
First, recall that $H$, as a closed subgroup, is a set of synthesis for $G$.
(\cite[Theorem 3]{TT2}).
\par
Suppose that there is $\sigma \in A(G)$ with $\sigma |_H = 1$ and
$\hat{\Theta}(\sigma)(T) = T$. Let $\phi \in A(G)$ be constant on
$H$. Then $\phi(e) \sigma - \phi$ vanishes on $H$, so that, by Lemma 
\ref{supplem2} we have,
\[
  0 = \hat{\Theta}(\phi(e) \sigma - \phi)(T) = \phi(e) \hat{\Theta}(\sigma)(T)
  - \hat{\Theta}(\phi)(T) = \phi(e) T - \hat{\Theta}(\phi)(T).
\]
This proves (i).
\par
For (ii), just note that, if $\phi \in \Mcb(A(G))$ is constant on $H$,
then $\phi(e) - \phi$ vanishes on $H$. An application of Corollary
\ref{suppcor1} then yields (ii).
\end{proof}
\par
As we already noted, our notion of support coincides with the one from 
\cite{Eym} for operators in $\VN(G)$. We now compute the support of
multiplication operators:
\begin{example}
Let $G$ be a locally compact group, and let $f \in L^\infty(G)$. 
For any $\phi \in A(G)$, we have
\begin{equation} \label{Linfty}
  \hat{\Theta}(\phi)(M_f) = M_f \hat{\Theta}(\phi)(1) = \phi(e) M_f ,
\end{equation}
where the first equality is due to the fact that $\hat{\Theta}(\phi)$ is an $L^\infty(G)$-bimodule homomorphism.
If $x \in G \setminus \{ e \}$, we can find $\phi \in A(G)$ with $\phi(x) 
\neq 0 = \phi(e)$, so that $\hat{\Theta}(\phi)(M_f) = 0$ by (\ref{Linfty}). 
Hence, $x$ cannot lie in $\supp \, M_f$. Since $x \in G \setminus \{ e \}$
was arbitrary, this means that $M_f \in {\cal B}_e(L^2(G))$. Moreover,
it is also clear from (\ref{Linfty}) that $e \in \supp \, M_f$ whenever
$f \neq 0$, i.e., $\supp \, M_f = \{ e \}$.
\end{example}
\section{Harmonic operators: the dual picture} \label{dualpic} 
Let $G$ be a locally compact group, and let $\Theta \!: M(G) \to
\CB({\cal B}(L^2(G)))$ be the completely isometric representation from 
\cite{Gha} (see also \cite{Neu}, \cite{Neu1}, and---for abelian 
$G$---\cite{Sto}) given by 
\[
  \langle \Theta (\mu) T , \omega \rangle := \int_G \langle \rho(t) 
   T \rho(t^{-1}) , \omega \rangle \, d\mu (t) 
  \qquad (\mu \in M(G), \,  T \in {\cal B}(L^2(G)), \, 
   \omega \in {\cal T}(L^2(G))) .
\]
For $\mu \in M^1(G)$, W.\ Jaworski and the first-named author define 
$T \in {\cal B}(L^2(G))$ to be \emph{$\mu$-harmonic} if $\Theta(\mu)(T) = T$ 
(\cite{JN}). Since
\[
  \Theta(\mu)M_\phi = M_{\mu \ast \phi} \qquad (\phi \in L^\infty(G)),
\]
this generalizes the notion of a $\mu$-harmonic function. More precisely, 
denoting the space of $\mu$-harmonic functions and operators 
by ${\cal H}_\mu$ and $\tilde{\cal H}_\mu$, respectively, 
this shows that ${\cal H}_\mu \subseteq \tilde{\cal H}_\mu$. Moreover, 
$\tilde{\cal H}_\mu$ carries a natural multiplication extending the one of 
${\cal H}_\mu$---and turning $\tilde{\cal H}_\mu$ into a non-commutative von 
Neumann algebra (\cite{JN}): given $S,T \in \tilde{\cal H}_\mu$, their product 
in $\tilde{\cal H}_\mu$ is explicitly given by 
\[
  \text{$w^\ast$-}\lim_n \int_G \rho(x) (ST) \rho (x^{-1}) \,\mu^{\ast n} (x). 
\]
One of the main results obtained in \cite{JN}---affirmatively answering a 
question by M.\ Izumi (\cite{Izu})---consists of a precise structural 
description of $\tilde{\cal H}_\mu$: provided 
that $G$ is second countable, $\tilde{\cal H}_\mu$ is exactly the crossed 
product of ${\cal H}_\mu$ with $G$ acting by left translation 
(\cite[Theorem 6.3]{JN}). In particular, it shows that the algebra $\VN(G)$, 
as a subalgebra of $\tilde{\cal H}_\mu$, plays the same r\^ole as the scalars 
for the classical algebra of harmonic functions ${\cal H}_\mu$. 
\par
Using $\hat{\Theta} \!: \Mcb(A(G)) \to \CB({\cal B}(L^2(G)))$ from
\cite{NRS}, we can extend the notion of a $\sigma$-harmonic functional
on $A(G)$ from \cite{CL} to that of a $\sigma$-harmonic operator in
a way analogous to the passage from ${\cal H}_\mu$ to $\tilde{\cal H}_\mu$ via $\Theta$:
\begin{definition}
Let $G$ be a locally compact group, and let $\sigma \in \Mcb(A(G))$.
Then an operator $T \in {\cal B}(L^2(G))$ is called
\emph{$\sigma$-harmonic} if $\hat{\Theta}(\sigma)(T) = T$. We denote
the collection of all $\sigma$-harmonic operators on $L^2(G)$ by
$\tilde{\cal H}_\sigma$.
\end{definition}
\begin{remarks}
\item Obviously, $\tilde{\cal H}_\sigma$ is a $w^\ast$-closed subspace
of ${\cal B}(L^2(G))$.
\item Trivially, $\tilde{\cal H}_\sigma$ contains ${\cal H}_\sigma$, and since,
\[
  \hat{\Theta}(\phi)(M_f) = M_f \, \hat{\Theta}(\phi)(1) = M_f \qquad (f \in L^\infty(G)),
\]
it contains $L^\infty(G)$ as well.
\end{remarks}
\par
As proven in \cite{CL}, ${\cal H}_\sigma$ is a von Neumann subalgebra
of $\VN(G)$ for $\sigma \in P^1(G)$, which stands in marked contrast to
${\cal H}_\mu$ with $\mu \in M^1(G)$. In the remainder of this section, 
we shall see that a similar statement is true for $\tilde{\cal H}_\sigma$:
under some, fairly mild, additional hypotheses, it is a von Neumann 
subalgebra of ${\cal B}(L^2(G))$; in fact, we shall prove that
\begin{equation} \label{generated}
  \tilde{\cal H}_\sigma = ({\cal H}_\sigma \cup L^\infty(G))''.
\end{equation}
\begin{remark}
The description (\ref{generated}) may be
viewed as a dual version of the central structural result 
\cite[Theorem 6.3]{JN}; in our setting, it is the algebra $L^\infty(G)$ which 
plays the r\^ole of the scalars. Noting that $L^\infty(G)$ and $\VN(G)$ are
Kac algebras dual to each other (see \cite{ES}), we are 
inclined to believe that one may find one single structure result that
unifies those descriptions of $\tilde{\cal H}_\sigma$ and $\tilde{\cal H}_\mu$ 
in the general framework of Kac algebras (or, even more generally, 
locally compact quantum groups). 
\end{remark}
\par 
We proceed through a series of lemmas and propositions.
\begin{proposition} \label{prop1}
Let $G$ be a locally compact group, and let $\sigma \in P^1(G)$. Then
\[
  ({\cal H}_\sigma \cup L^\infty(G))'' \subset \tilde{\cal H}_\sigma
\]
holds.
\end{proposition}
\begin{proof}
By von Neumann's double commutant theorem, $({\cal H}_\sigma \cup L^\infty(G))''$ is the von
Neumann subalgebra of ${\cal B}(L^2(G))$ generated by ${\cal H}_\sigma$
and $L^\infty(G)$. Since $\tilde{\cal H}_\sigma$ is a $w^\ast$-closed
subspace of ${\cal B}(L^2(G))$, it is sufficient to show that finite
products of operators from ${\cal H}_\sigma \cup L^\infty(G)$ belong to
$\tilde{\cal H}_\sigma$. Moreover, since ${\cal H}_\sigma \cong \VN(G_\sigma)$, it is enough to consider products of the form
\begin{equation} \label{product}
  \prod_{j=1}^n \lambda(x_j) M_{f_j},
\end{equation}
where $x_1, \ldots, x_n \in G_\sigma$ and $f_1, \ldots, f_n \in L^\infty(G)$. Since
\[
  \lambda(x) M_f = M_{\lambda(x) f} \lambda(x) \qquad (x \in G, \, f \in L^\infty(G)),
\]
a simple induction on the number of factors in (\ref{product}) shows that (\ref{product}) is, in fact, of the form $M_f \lambda(x)$ with $f \in L^\infty(G)$ and $x \in G_\sigma$. Since $\hat{\Theta}(\sigma)$ is an $L^\infty(G)$-bimodule homomorphism and
since $\hat{\Theta}(\sigma)$ fixes $\lambda(x)$ if $x \in G_\sigma$, we have
\[
  \hat{\Theta}(\sigma)(M_f \lambda(x)) = 
  M_f \hat{\Theta}(\sigma)(\lambda(x)) =  M_f \lambda(x)  
  \qquad (f \in L^\infty(G), \, x \in G_\sigma).
\]
i.e., $M_f \lambda(x) \in \tilde{\cal H}_\sigma$.
\par
In view of the foregoing remarks, this completes the proof.
\end{proof}
\par
While Proposition \ref{prop1} provides an estimate for $\tilde{\cal H}_\sigma$ ``from below'', we now give one ``from above''. 
\par 
Extending our previous notation, we set $G_\sigma := \{ x \in G : \sigma(x) = 1 \}$ for \emph{any} $\sigma \in \Mcb(A(G))$. Note that
$G_\sigma$ need not be a subgroup of $G$ unless $\sigma \in P^1(G)$, but that it is still a closed subset of $G$, so that the conclusion of the following proposition still makes sense.
\begin{proposition} \label{prop2}
Let $G$ be a locally compact group, and let $\sigma \in \Mcb(A(G))$. Then
we have $\tilde{\cal H}_\sigma \subset {\cal B}_{G_\sigma}(L^2(G))$.
\end{proposition}
\begin{proof}
Let $T \in \tilde{\cal H}_\sigma$, and let $x \in G \setminus G_\sigma$, i.e.,
$\sigma(x) \neq 1$. Let $\phi \in A(G)$ be such that $\phi(x) \neq 0$, so that
$\psi(x) := \phi(x)\sigma(x) - \phi(x) \neq 0$. On the other hand, we have
\[
  \hat{\Theta}(\psi)(T) = 
  \hat{\Theta}(\phi\sigma)(T) - \hat{\Theta}(\phi)(T) = 
  \hat{\Theta}(\phi) (\hat{\Theta}(\sigma)(T)) - \hat{\Theta}(\phi)(T) 
  = \hat{\Theta}(\phi)(T) - \hat{\Theta}(\phi)(T) =0.
\]
This means that $x \notin \supp \, T$.
\end{proof}
\par
Let $G$ be a locally compact group, and let $H$ be a closed subgroup. 
We define
\[
  L^\infty(G:H) := \{ \phi \in L^\infty(G) : 
                   \text{$\lambda(x)\phi = \phi$ for all $x \in H$} \}.
\]
It is obvious that $L^\infty(G:H)$ is $w^\ast$-closed in $L^\infty(G)$. 
Slightly deviating from \cite{For}, we set
\[
  B(G:H) := \{ \phi \in B(G) : 
                   \text{$\lambda(x)\phi = \phi$ for all $x \in H$} \}.
\]
(Note that in \cite{For}, right instead of left translates are considered.)
We can canonically embed $B(G:H)$ into $L^\infty(G:H)$. Also, we set
\[
  P^1_H(G) := \{ \sigma \in P^1(G) : \sigma |_H \equiv 1 \}.
\]
As is observed in \cite{KL}, the functions in
$P^1_H(G)$ are constant on both left and right cosets of $H$ and thus
are contained in $B(G:H)$.
Following \cite{KL}, we say that $G$ has the \emph{$H$-separation property} 
if, for each $x \in G \setminus H$, there is $\sigma \in P^1_H(G)$ such that
$\sigma(x) \neq 1$. For instance, whenever $H$ is open, compact, or 
neutral---this includes all normal subgroups as well as all closed subgroups 
of $[\SIN]$-groups---, $G$ has the $H$-separation property (see \cite{KL}).
\begin{lemma} \label{wdense}
Let $G$ be a locally compact group, and let $H$ be a closed subgroup of $G$
such that $G$ has the $H$-separation property.
Then $B(G:H)$ is $w^\ast$-dense in $L^\infty(G:H)$.
\end{lemma}
\begin{proof}
Let $\M$ denote the $w^\ast$-closure of the linear span of all
right translates of all functions in $P^1_H(G)$ in $L^\infty(G)$.
Then $\M$ is a von Neumann subalgebra of $L^\infty(G)$ and invariant 
under right translation. Set
\[
  L := \{ x \in G : \text{$\lambda(x)\phi = \phi$ for all $\phi \in \M$} \}.
\]
Then $L$ is a closed subgroup of $G$ containing $H$, and from 
\cite[Theorem 2]{TT} (see also \cite[Lemma 3.2]{Lau}), we conclude that
$\M = L^\infty(G:L)$.
\par
Assume that there is $x \in L \setminus H$. Since $G$ has the $H$-separation
property, there is $\sigma \in P^1_H(G) \subset L^\infty(G:L)$ such that 
$\sigma(x) \neq 1$ and thus $\lambda(x^{-1}) \sigma \neq \sigma$. This is
a contradiction, so that $H = L$.
\end{proof}
\par
For any locally compact group $G$ and $\sigma \in P^1(G)$, it is clear by
definition that $G$ has the $G_\sigma$-separation property. Hence, we obtain:
\begin{corollary} \label{wdensecor} 
Let $G$ be a locally compact group, and let $\sigma \in P^1(G)$.
Then $B(G:G_\sigma)$ is $w^\ast$-dense in $L^\infty(G:G_\sigma)$.
\end{corollary}
\begin{lemma} \label{invlem}
Let $G$ be a locally compact group, let $H$ be a closed subgroup, and let
$\phi \in B(G:H)$. Then there are a unitary representation $\pi$ of $G$
on a Hilbert space $\Hilbert$ as well as $\xi, \eta \in \Hilbert$ such that
\[
  \phi(x) = \langle \pi(x) \xi, \eta \rangle \qquad (x \in G)
\]
and $\pi(y)\eta = \eta$ for all $y \in H$.
\end{lemma}
\begin{proof}
By definition of $B(G)$, there are a unitary representation $\pi$ of
$G$ on a Hilbert space $\Hilbert$ as well as $\xi, \eta_0 \in \Hilbert$
such that
\[
  \phi_{\xi,\eta_0}(x) := \langle \pi(x) \xi, \eta_0 \rangle 
  = \phi(x) \qquad (x \in G).
\]
Without loss of generality, suppose that $\| \xi \| = 1$. Set
\[
  C := \{ \eta \in \Hilbert : \text{$\| \eta \| \leq \| \eta_0 \|$
  and $\phi_{\xi,\eta} = \phi$} \}.
\]
Then $C$ is a non-empty, convex, weakly compact subset of $\Hilbert$.
Let $\eta \in C$ and $y \in H$. Then we have
\[
  \phi_{\xi,\pi(y)\eta}(x) = \langle \pi(x) \xi, \pi(y)\eta \rangle \\
  = \langle \pi(y^{-1} x) \xi, \eta \rangle = (\lambda(y)\phi)(x)
  = \phi(x) \qquad (x \in G),
\]
so that $\pi(y)\eta \in C$ again. From the Ryll-Nardzewski 
fixed point theorem (\cite[(9.6) Theorem]{GD}), we 
conclude that there is $\eta \in C$ with $\pi(y)\eta = \eta$ for $y \in H$. 
\end{proof}
\begin{proposition} \label{supprop}
Let $G$ be a locally compact group, let $H$ be a closed subgroup of 
$G$, let $T \in {\cal B}_H(L^2(G))$, and suppose that one of the following 
holds:
\begin{alphitems}
\item there is $\sigma \in A(G)$ with $\sigma |_H = 1$ and 
$\hat{\Theta}(\sigma)(T) = T$;
\item $G$ has the approximation property and the $H$-separation property.
\end{alphitems}
Then $T$ lies in $(\VN(H) \cup L^\infty(G))''$.
\end{proposition}
\par
For the proof, we first establish some conventions for (possibly infinite) matrices. The extensive use of matrix calculations may seem like an unnecessary complication at the first glance, but appears to be unavoidable in view of how $\hat{\Theta} \!: \Mcb(A(G)) \to \CB({\cal B}(L^2(G))$ is defined in \cite{NRS}, namely via the extended Haagerup tensor product of operator spaces, an object that itself is defined in terms of arbitrarily large matrices.
\par 
Given a linear space $E$ and index sets $\mathbb I$ and $\mathbb J$,
we denote by $M_{\mathbb{I} \times \mathbb{J}}(E)$ the matrices $[ x_{i,j}
]_{i \in \mathbb{I} \atop j \in \mathbb J}$ with $x_{i,j} \in E$ for
$(i,j) \in \mathbb{I} \times \mathbb{J}$. If $\mathbb{I} = \mathbb J$, we
write $M_{\mathbb J}(E)$, and if $E = \comps$, we simply use the symbols
$M_{\mathbb{I} \times \mathbb{J}}$ and $M_{\mathbb J}$ instead of
$M_{\mathbb{I} \times \mathbb{J}}(\comps)$ and $M_{\mathbb J}(\comps)$,
respectively. For $x \in E$ and $[ \alpha_{i,j} 
]_{i \in \mathbb{I} \atop j \in \mathbb J}$, we set
$[ \alpha_{i,j} ]_{i \in \mathbb{I} \atop j \in \mathbb J} \tensor x =
[ \alpha_{i,j} x ]_{i \in \mathbb{I} \atop j \in \mathbb J}$. We convene
to interpret families $[ x_j ]_{j \in \mathbb J}$ as row vectors, and write 
$[ x_j ]_{j \in \mathbb J}^t$ for the corresponding column vector.
Finally, we denote by $1_\mathbb{J}$ the matrix $[ \delta_{j,k} ]_{j,k \in
\mathbb J}$, and set $0_{\mathbb{I} \times \mathbb{J}} = 
[ \alpha_{i,j}]_{i \in \mathbb{I} \atop j \in \mathbb J}$ with $\alpha_{i,j}
= 0$ for $(i,j) \in \mathbb{I} \times \mathbb{J}$.
\begin{proof}
If (a) holds, $H$ is necessarily compact, so that in both case (a) and 
case (b), $G$ has the $H$-separation property. By Lemma
\ref{wdense}, $B(G:H)$ is therefore $w^\ast$-dense in $L^\infty(G:H)$.
We will show that $M_\phi T = T M_\phi$ for
all $\phi \in B(G:H)$, so that $T$ lies in $L^\infty(G:H)'$. Since,
as is easily checked,
\[
  (\VN(H) \cup L^\infty(G))' = L^\infty(G:H),
\]
this will prove the proposition. 
\par
We adapt part of the proof of \cite[Theorem 5.1]{NRS} to our situation.
\par
Let $\phi \in B(G:H)$. By Lemma \ref{invlem}, there thus are a unitary 
representation $\pi$ of $G$ on some Hilbert space $\Hilbert$ as well as 
$\xi, \eta \in \Hilbert$ such that
\[
  \phi(x) = \langle \pi(x) \xi, \eta \rangle \quad (x \in G)
\]
and $\pi(y) \eta = \eta$ for $y \in H$. Let 
\[
  \mathfrak{K} := \{ \zeta \in \Hilbert : \text{$\pi(y) \zeta = \zeta$
  for $y \in H$} \},
\]
so that, in particular, $\eta \in \mathfrak K$.
Let $( e_i )_{i \in \mathbb I}$ be an orthonormal basis for $\mathfrak K$,
and extend it to an orthonormal basis $( e_i )_{i \in \mathbb J}$ for $\Hilbert$ (so that, in particular, $\mathbb{I} \subset \mathbb{J}$). For $(i,k) \in \mathbb{J} \times \mathbb J$ set
\[
  \phi_{i,k}(x) := \langle \pi(x) e_k, e_i \rangle 
  \quad\text{and}\quad
  \check{\phi}_{i,k}(x) := \phi_{i,k}(x^{-1}) 
  \qquad (x \in G)   
\]
Since $e_i \in \mathfrak K$, it follows that $\phi_{i,k} \in B(G:H)$ for all $(i,k) \in \mathbb{I} \times \mathbb J$. By the definition
of $\hat{\Theta}$, we have
\begin{equation} \label{fart1}
  \hat{\Theta}(\phi_{i,k})(T) = 
  \sum_{j \in \mathbb J} M_{\phi_{i,j}} T M_{\check{\phi}_{j,k}} 
  \qquad ((i,k) \in \mathbb{I} \times \mathbb{J});
\end{equation}
By Proposition \ref{supprop1}(i) or (ii)---depending on whether (a) or (b) is satisfied---, the left hand side of (\ref{fart1})
equals $\phi_{i,k}(e)T$ for all $(i,k) \in \mathbb{I} \times \mathbb{J}$, so that (\ref{fart1}) becomes
\begin{equation} \label{fart+}
  \phi_{i,k}(e)T = \sum_{j \in \mathbb J} M_{\phi_{i,j}} T M_{\check{\phi}_{j,k}} 
  \qquad ((i,k) \in \mathbb{I} \times \mathbb{J}).
\end{equation}
Interpreting (\ref{fart+}) as a matrix identity---with matrices indexed
over $\mathbb{J} \times \mathbb{J}$---, we obtain
\begin{equation} \label{fart2}
  \phi_{i,k}(e)T  = [ \phi_{i,j}]_{j \in \mathbb J} \,
  (1_\mathbb{J} \tensor T) \,   \left[ \check{\phi}_{j,k} \right]_{j \in \mathbb J}
  \qquad ((i,k) \in \mathbb{I} \times \mathbb{J}).
\end{equation}
We view $\Phi := [ \phi_{i,k} ]_{i,k \in \mathbb J}$ as an element of 
$M_\mathbb{J}({\cal C}_b(G))$ and set $\check{\Phi} := 
\left[ \check{\phi}_{i,k} \right]_{i,k \in \mathbb J}$. Then 
$\check{\Phi}$ also
lies in $M_\mathbb{J}({\cal C}_b(G))$ and satisfies
\begin{equation} \label{unitary}
  \Phi \, \check{\Phi} = 1_\mathbb{J} \tensor 1 = \check{\Phi} \, \Phi.
\end{equation}
Furthermore, set $\Psi := 
[ \phi_{i,k} ]_{i \in \mathbb{I} \atop k \in \mathbb J} \in
M_{\mathbb{I} \times \mathbb J}({\cal C}_b(G))$.
Since $\phi_{i,k}(e) = \delta_{i,k}$ for $(i,k)\in \mathbb{I}
\times \mathbb J$, we obtain 
from (\ref{fart2}) that
\[
  \begin{split}
  \left(1_\mathbb{I} \oplus 
  0_{\mathbb{I} \times (\mathbb{J} \setminus \mathbb{I})} \right)
  \tensor T & 
  = [ \phi_{i,k}(e)T ]_{i \in \mathbb{I} \atop k \in \mathbb J} \\
  & = \left[ [ \phi_{i,j}]_{j \in \mathbb J} \,
  (1_\mathbb{J} \tensor T) \, 
  \left[ \check{\phi}_{j,k} \right]_{j \in \mathbb J} 
  \right]_{i \in \mathbb{I} \atop k \in \mathbb J} \\ 
  & = \Psi \, ( 1_\mathbb{J} \tensor T) \, \check{\Phi},
  \end{split}
\]
and thus, by (\ref{unitary})
\begin{equation} \label{commutes}
  \Psi \, ( 1_\mathbb{J} \tensor T) = 
  \left( \left(1_\mathbb{I} \oplus 
  0_{\mathbb{I} \times (\mathbb{J} \setminus \mathbb{I})}\right)
  \tensor T \right) \, \Phi.
\end{equation}
Let $[ \alpha_j ]_{j \in \mathbb J}$ and $[ \beta_i ]_{i \in \mathbb I}$ be
in $\ell^2(\mathbb{J})$ and $\ell^2(\mathbb{I})$, respectively, such that 
$\xi = \sum_{j \in \mathbb J} \alpha_j e_j$
and $\eta = \sum_{i \in \mathbb I} \beta_i e_i$; it follows that
\[
  \phi(x) = \langle \pi(x)\xi, \eta \rangle =
  \sum_{j \in \mathbb{J} \atop i \in \mathbb I} 
  \alpha_j \langle \pi(x)e_j,e_i \rangle \bar{\beta}_i=
  \sum_{j \in \mathbb{J} \atop i \in \mathbb I} 
  \alpha_j \, \phi_{i,j}(x) \, \bar{\beta}_i
  \qquad (x \in G)
\]
or, in matrix notation,
\[
  \phi = \left[ \bar{\beta}_i \right]^t_{i \in \mathbb I}
  \, \Psi \, [ \alpha_j ]_{j \in \mathbb J}.
\]
Hence, we obtain eventually:
\[
  \begin{split}
  M_\phi T & =  \left[ \bar{\beta}_i \right]_{i \in \mathbb I}^t \, \Psi \, 
  (1_\mathbb{J} \tensor T) \, [ \alpha_j ]_{j \in \mathbb J} \\
  & = \left[ \bar{\beta}_i \right]_{i \in \mathbb I}^t \,
  \left( \left(1_\mathbb{I} \oplus 
  0_{\mathbb{I} \times (\mathbb{J} \setminus \mathbb{I})} \right)
  \tensor T \right) \, \Phi
  \,  [ \alpha_j ]_{j \in \mathbb J} \\
  & = TM_\phi.
  \end{split}
\]
This proves the claim.
\end{proof}
\par
We can now prove the main result of this section:
\begin{theorem}
Let $G$ be a locally compact group, and let $\sigma \in P^1(G)$. 
Then the following hold:
\begin{items}
\item if $\sigma \in A(G)$, then $\tilde{\cal H}_\sigma =
({\cal H}_\sigma \cup L^\infty(G))''$;
\item if $G$ has the approximation property, then
\[
  \tilde{\cal H}_\sigma = ({\cal H}_\sigma \cup L^\infty(G))'' 
  = {\cal B}_{G_\sigma}(L^2(G)).
\]
\end{items}
In either case, $\tilde{\cal H}_\sigma =  ({\cal H}_\sigma \cup L^\infty(G))''$ 
is a von Neumann subalgebra of ${\cal B}(L^2(G))$.
\end{theorem}
\begin{proof}
By Propositions \ref{prop1} and \ref{prop2} 
\[
  ({\cal H}_\sigma \cup L^\infty(G))'' \subset \tilde{\cal H}_\sigma 
  \subset {\cal B}_{G_\sigma}(L^2(G))
\]
holds without any additional hypotheses.
\par
In case (i), we conclude from Proposition \ref{supprop} (with condition
(a)) that $\tilde{\cal H}_\sigma \subset 
({\cal H}_\sigma \cup L^\infty(G))''$. 
For case (ii), recall that, as we remarked before Corollary 
\ref{wdensecor}, $G$ has the $G_\sigma$-separation property.
Hence, by Proposition \ref{supprop}  (with condition (b)), we even 
have ${\cal B}_{G_\sigma}(L^2(G)) 
\subset ({\cal H}_\sigma \cup L^\infty(G))''$. This proves the theorem.
\end{proof}
\begin{remark}
We believe, but have been unable to prove, that
\[
  \tilde{\cal H}_\sigma = ({\cal H}_\sigma \cup L^\infty(G))'' 
  = {\cal B}_{G_\sigma}(L^2(G))
\]
holds without any additional hypotheses on $\sigma$ or $G$.
\end{remark}
\par
For adapted $\sigma$, we obtain as a special case:
\begin{corollary}
Let $G$ be a locally compact group, and let $\sigma \in P^1(G)$ be
adapted. Then the following hold:
\begin{items}
\item if $\sigma \in A(G)$, then $\tilde{\cal H}_\sigma =
L^\infty(G)$;
\item if $G$ has the approximation property, then
\[
  \tilde{\cal H}_\sigma = L^\infty(G) = {\cal B}_e(L^2(G)).
\]
\end{items}
In either case, $\tilde{\cal H}_\sigma = L^\infty(G)$ is a von
Neumann subalgebra of ${\cal B}(L^2(G))$.
\end{corollary} 
\section{Harmonic operators through ideals in the predual} \label{ideals} 
An alternative, very fruitful approach to the classical space ${\cal H}_\mu$ of 
$\mu$-harmonic functions (where $\mu$ is a probability measure) has been carried out 
by G.\ A.\ Willis (\cite{Wil}). There, the study of ${\cal H}_\mu$ is 
transferred to its pre-annihilator 
\[
  J_\mu := \overline{\{ f - f * \mu : f \in L^1(G) \}}^{L^1(G)}
\]
in $L^1(G)$, which (obviously) forms a left ideal in the group algebra. In 
order to investigate the non-commutative analogue 
$\tilde{\cal H}_\mu$ of ${\cal H}_\mu$ in a 
similar fashion, W.\ Jaworski and the first-named author use the ``quantized'' 
convolution product in the space ${\cal T}(L^2(G)) = {\cal B}(L^2(G))_\ast$ 
of trace class operators, as introduced and studied in \cite{Neu} and 
\cite{Neu2} (see also \cite{Pir} for further results on this product). Indeed, 
they show that the pre-annihilator 
\[
  \tilde{J}_\mu := \overline{\{ \omega - \Theta(\mu)_\ast (\omega) : 
  \omega \in {\cal T}(L^2(G)) \}}^{{\cal T}(L^2(G))}
\]
of $\tilde{\cal H}_\mu$ in ${\cal T}(L^2(G))$ forms a left ideal with respect 
to this convolution (\cite[Proposition 3.3]{JN}). 
\par 
We shall first indicate how to equip ${\cal T}(L^2(G))$ with a product 
that turns it into a completely contractive Banach algebra, and may be thought 
of as a ``non-commutative pointwise product'', in other words, a Fourier 
algebra type product, instead of a convolution type product, as described 
above. Just as the latter, the multiplication we shall consider here is also 
very different from the ordinary composition of operators. From a 
Hopf--von Neumann algebraic point of view, both the convolution type and the 
Fourier algebra type product on the space ${\cal T}(L^2(G))$ are constructed 
in exactly the same way: the first one is derived from a co-multiplication 
on ${\cal B}(L^2(G))$ that canonically extends the one of $L^\infty(G)$
whereas the second one is based on a canonical extension of the 
co-multiplication of $\VN(G)$. In this sense, 
both products on ${\cal T}(L^2(G))$ are ``dual'' to each other. 
\par 
Let $G$ be a locally compact group. We define a unitary operator
$W \in {\cal B}(L^2(G \times G))$ by letting
\[
  (W \xi)(x,y) := \xi(x,xy) \qquad (\xi \in L^2(G), \, x,y \in G).
\]
Identifying $L^2(G \times G)$ with $L^2(G) \ttensor_2 L^2(G)$ (Hilbert space
tensor product), we denote the flip map on $L^2(G \times G)$ by $\sigma$.
Then $\hat{W} := \sigma W^\ast \sigma$ is a again unitary, and
\[
  \hat{\Gamma} \!: {\cal B}(L^2(G)) \to {\cal B}(L^2(G \times G)),
  \quad T \mapsto \hat{W}(1 \tensor T) \hat{W}^\ast
\]
is a \emph{co-multiplication}, i.e., a normal, unital, injective 
$^\ast$-homomorphism satisfying
\[
  (\hat{\Gamma} \tensor \id) \circ \hat{\Gamma} = 
  (\id \tensor \hat{\Gamma}) \circ \hat{\Gamma}. 
\]
The pre-adjoint
\[
  \hat{\Gamma} \!: {\cal T}(L^2(G)) \Tensor {\cal T}(L^2(G)) \to
  {\cal T}(L^2(G)),
\]
where $\Tensor$ denotes the projective tensor product of operator spaces,
is a complete contraction and turns ${\cal T}(L^2(G))$ into a 
completely contractive 
Banach algebra. In order to tell this product on ${\cal T}(L^2(G))$
apart from the usual composition of operators, we denote it by $\bullet$,
so that
\[
  \langle \omega \bullet \rho, T \rangle := \langle \omega \tensor \rho,
  \hat{\Gamma}(T) \rangle \qquad (\omega, \rho \in {\cal T}(L^2(G)), \,
  T \in {\cal B}(L^2(G))).
\]
Furthermore, the pre-adjoint of the inclusion $\VN(G) \subset {\cal B}(L^2(G))$
is an algebra homomorphism $\pi \!: {\cal T}(L^2(G)) \to A(G)$ 
(and necessarily a complete quotient map), and ${\cal B}(L^2(G))$ becomes a completely
contractive ${\cal T}(L^2(G))$-bi\-mo\-dule through
\begin{multline*}
  \omega \cdot T := (\id \tensor \omega)(\hat{\Gamma}(T))
  \quad\text{and}\quad
  T \cdot \omega := (\omega \tensor \id)(\hat{\Gamma}(T)) \\
  \qquad (\omega \in {\cal T}(L^2(G)), \, T \in {\cal B}(L^2(G))).
\end{multline*}
(These module actions have been studied in \cite{NRS} and \cite{PR}.) 
\par 
Since $\hat{W} \in \VN(G) \bar{\tensor} L^\infty(G)$, we have 
\[
  \omega \cdot T = (\id \tensor \omega)(\hat{W}(1 \tensor T)\hat{W}^\ast) \in \VN(G)
  \qquad (\omega \in {\cal T}(L^2(G)), \, T \in {\cal B}(L^2(G))). 
\]
Hence, the module action induces a complete contraction $\hat{\mathcal{S}} \!: {\cal T}(L^2(G)) 
\Tensor {\cal B}(L^2(G)) \to \VN(G)$. In \cite{PR}, 
$\mathrm{UCQ}(G)$ is defined to be the range of $\hat{\mathcal{S}}$ in 
$\VN(G)$, equipped with the quotient operator space structure. 
The dual space $\mathrm{UCQ}(G)^\ast$ then naturally
becomes a completely contractive Banach algebra, which
contains $\Mcb A(G)$ as a closed subalgebra, and the adjoint 
$\hat{\mathcal{S}}^\ast$ extends 
the representation $\hat{\Theta}$ from $\Mcb(A(G))$ to $\mathrm{UCQ}(G)^\ast$. 
More precisely, if, for $n \in \mathrm{UCQ}(G)^\ast$ and 
$T \in {\cal B}(L^2(G))$, one defines $n \diamond T \in {\cal B}(L^2(G))$ 
through 
\[
  \langle n \diamond T , \omega \rangle := 
  \langle n , \omega \cdot T \rangle \qquad 
  (\omega \in {\cal T}(L^2(G))) 
\] 
then, for each $n \in \mathrm{UCQ}(G)^\ast$, the map 
\[
  {\cal B}(L^2(G)) \to {\cal B}(L^2(G)), \quad T \mapsto n \diamond T
\] 
is a completely bounded operator on ${\cal B}(L^2(G))$, which we denote by 
$\tilde{\hat{\Theta}}(n)$. It is easy to check that $\tilde{\hat{\Theta}} = 
\hat{\mathcal{S}}^\ast$. For more information, see \cite{PR} 
(and \cite{NRS} for the amenable case). 
\par 
In the sequel, we shall use the fact that $\tilde{\hat{\Theta}}$ and 
$\hat{\Theta}$ coincide on $A(G)$, which follows from \cite{PR} 
(or \cite{NRS} if $G$ is amenable). For the reader's convenience, we include 
a different (and short) proof:
\begin{lemma} \label{thetatheta} 
Let $G$ be a locally compact group. Then we have 
\[
  \tilde{\hat{\Theta}} (\phi) = \hat{\Theta} (\phi) \qquad (\phi \in A(G)).
\]
\end{lemma}
\begin{proof} 
Let $\phi \in A(G)$, and note that both $\tilde{\hat{\Theta}} (\phi)$ and 
$\hat{\Theta} (\phi)$ are normal, which follows easily from the definition 
(see also \cite[Theorem 4.3]{NRS} for $\hat{\Theta} (\phi)$). 
\par 
Both $\tilde{\hat{\Theta}} (\phi)$ and $\hat{\Theta} (\phi)$ are 
$L^\infty(G)$-bimodule maps: we already know this for $\hat{\Theta} (\phi)$, and it follows for $\tilde{\hat{\Theta}} (\phi)$ 
from \cite[Theorem 2.3]{PR}. Furthermore,  $\tilde{\hat{\Theta}} (\phi)$ and $\hat{\Theta} (\phi)$ coincide on all operators 
$\lambda(x)$ with $x \in G$: we know that $\hat{\Theta} (\phi)$ on $\VN(G)$ 
is nothing but the canonical action of $\phi$ on $\VN(G)$, and the same follows for $\tilde{\hat{\Theta}} (\phi)$ from its definition.
Consequently, we have
\begin{multline*}
  \tilde{\hat{\Theta}} (\phi)(M_f \lambda(x)) = M_f \tilde{\hat{\Theta}} (\phi)(\lambda(x)) \\
  = M_f \hat{\Theta} (\phi)(\lambda(x)) = \hat{\Theta} (\phi)(M_f \lambda(x)) \qquad (f \in L^\infty(G), \, x \in G).
\end{multline*}
\par
As in the proof of Proposition \ref{prop1}, we conclude that $\tilde{\hat{\Theta}} (\phi) = \hat{\Theta} (\phi)$.
\end{proof} 
Following \cite{CL}, where the analogous questions for $\VN(G)$ and
${\cal H}_\sigma$ were considered, we now study the pre-annihilator
of $\tilde{\cal H}_\sigma$ in ${\cal T}(L^2(G))$ 
for any $\sigma \in \Mcb(A(G))$.
\begin{lemma} \label{rightmod}
For any locally compact group $G$, we have
\[
  T \cdot \omega = \hat{\Theta}(\pi(\omega))(T) 
  \qquad (\omega \in {\cal T}(L^2(G)), \, T \in {\cal B}(L^2(G)).
\]
where $\pi \!: {\cal T}(L^2(G)) \to A(G)$ is the canonical quotient map. 
\end{lemma}
\begin{proof}
Let $\omega \in {\cal T}(L^2(G))$, and let $T \in {\cal B}(L^2(G))$.
Then we have
\begin{equation*} \label{modid}
  \langle \rho, T \cdot \omega \rangle = \langle \omega\bullet\rho, T \rangle 
  = \langle \omega, \rho \cdot T \rangle
  \qquad (\rho \in {\cal T}(L^2(G))).
\end{equation*}
Since $\rho \cdot T \in \VN(G)$ for
$\rho \in {\cal T}(L^2(G))$, the evaluation of $\omega$ at 
$\rho \cdot T$ depends only on $\pi(\omega)$,
so that $\langle \rho, T \cdot \omega \rangle = \langle \pi(\omega), \rho \cdot T \rangle$ 
for each $\rho \in {\cal T}(L^2(G))$. Hence we obtain 
\[ 
  \begin{split} 
  \langle \rho, T \cdot \omega \rangle 
  & = \langle \pi(\omega), \rho \cdot T \rangle \\ 
  & = \langle \pi(\omega), \hat{\mathcal{S}} (\rho \otimes T) \rangle \\ 
  & = \langle \rho \otimes T , \hat{\mathcal{S}}^* (\pi(\omega)) \rangle \\ 
  & = \langle \rho \otimes T , \tilde{\hat{\Theta}} (\pi(\omega)) \rangle \\ 
  & = \langle \rho , \tilde{\hat{\Theta}} (\pi(\omega)) (T) \rangle \\ 
  & = \langle \rho , \hat{\Theta} (\pi(\omega)) (T) \rangle, \qquad\text{by Lemma \ref{thetatheta}}, 
  \end{split} 
\] 
for all $\rho \in {\cal T}(L^2(G))$, as desired. 
\end{proof}
The following is a consequence of \cite[Lemma 2.1]{PR}, using the fact that 
$\tilde{\hat{\Theta}} \mid_{\Mcb((A(G))} = \hat{\Theta}$. We include our own (short) proof for the sake of completeness:
\begin{lemma} \label{leftmod}
For any locally compact group $G$, we have
\[
  \hat{\Theta}(\phi)(\omega \cdot T) = \omega \cdot \hat{\Theta}(\phi)(T)
  \qquad (\phi \in \Mcb((A(G)), \,
  \omega \in {\cal T}(L^2(G)), \, T \in {\cal B}(L^2(G))).
\]
\end{lemma}
\begin{proof}
Let $\phi \in \Mcb((A(G))$, let $\omega, \rho \in {\cal T}(L^2(G))$, 
and let $T \in {\cal B}(L^2(G))$. 
\par 
We first note that
\[
  \begin{split}
  \langle \pi(\hat{\Theta}(\phi)_\ast(\rho)), S \rangle 
  & = \langle \hat{\Theta}(\phi)_\ast(\rho), S \rangle \\
  & = \langle \rho, \hat{\Theta}(\phi)(S) \rangle \\
  & = \langle \rho, \phi \cdot S \rangle \\
  & = \langle \pi(\rho), \phi \cdot S \rangle \\
  & = \langle \pi(\rho) \phi, S \rangle \qquad (S \in \VN(G)),
  \end{split}
\]
so that
\begin{equation} \label{yai}
  \pi(\hat{\Theta}(\phi)_\ast(\rho)) = \pi(\rho) \phi.
\end{equation}
Now we obtain:
\[
  \begin{split}
  \langle \rho, \hat{\Theta}(\phi)(\omega \cdot T) \rangle
  & = \langle \hat{\Theta}(\phi)_\ast(\rho), \omega \cdot T \rangle \\
  & = \langle \hat{\Theta}(\phi)_\ast(\rho) \bullet \omega , T \rangle \\
  & = \langle \omega , T \cdot \hat{\Theta}(\phi)_\ast(\rho)\rangle \\
  & = \langle \omega, \hat{\Theta}(\pi(\hat{\Theta}(\phi)_\ast(\rho)))(T) 
      \rangle, \qquad\text{by Lemma \ref{rightmod}}, \\
  & = \langle \omega, \hat{\Theta}(\pi(\rho)\phi)(T) \rangle,
      \qquad\text{by (\ref{yai})}, \\
  & = \langle \omega, \hat{\Theta}(\pi(\rho))(\hat{\Theta}(\phi)(T)) \rangle \\
  & = \langle \omega, \hat{\Theta}(\phi)(T) \cdot \rho \rangle, \qquad
      \text{by Lemma \ref{rightmod} again}, \\
  & = \langle \rho \bullet \omega, \hat{\Theta}(\phi)(T) \rangle \\
  & = \langle \rho, \omega \cdot \hat{\Theta}(\phi)(T) \rangle.
  \end{split}
\]
Since $\rho \in {\cal T}(L^2(G))$ is arbitrary, this yields that
\[
  \hat{\Theta}(\phi)(\omega \cdot T) = \omega \cdot \hat{\Theta}(\phi)(T),
\]
as claimed.
\end{proof}
\begin{theorem}
Let $G$ be a locally compact group, and let $\sigma \in \Mcb(A(G))$.
Then
\[
  \tilde{I}_\sigma := \varcl{\{ \hat{\Theta}(\sigma)_\ast(\omega) - \omega :
  \omega \in {\cal T}(L^2(G)) \}}^{{\cal T}(L^2(G))}
\]
is the pre-annihilator of $\tilde{\cal H}_\sigma$ in ${\cal T}(L^2(G))$---so that $\tilde{\cal H}_\sigma \cong ({\cal T}(L^2(G)) / \tilde{I}_\sigma)^\ast$---and a two-sided ideal of $({\cal T}(L^2(G)), \bullet)$.
\end{theorem}
\begin{proof}
It is straightforward to see that 
\[
  \omega \in \tilde{I}_\sigma \iff \text{$\langle \omega, T \rangle = 0$ for all $T \in \tilde{\cal H}_\sigma$}
  \qquad (\omega \in {\cal T}(L^2(G))).
\]
Hence, $\tilde{I}_\sigma$ is indeed the pre-annihilator of $\tilde{\cal H}_\sigma$ in ${\cal T}(L^2(G))$, and $({\cal T}(L^2(G)) / \tilde{I}_\sigma)^\ast = \tilde{I}_\sigma^\perp \cong \tilde{\cal H}_\sigma$ holds.
\par
We first show that $\tilde{I}_\sigma$ is a left ideal in ${\cal T}(L^2(G))$.
To this end, let $\omega, \rho \in {\cal T}(L^2(G))$, and note that
\[
  \begin{split} 
  \langle \omega \bullet \hat{\Theta}(\sigma)_\ast(\rho), T \rangle & =
  \langle \hat{\Theta}(\sigma)_\ast(\rho), T \cdot \omega \rangle \\
  & = \langle \rho, \hat{\Theta}(\sigma)(T \cdot \omega) \rangle \\
  & = \langle \rho, \hat{\Theta}(\sigma)(\hat{\Theta}(\pi(\omega))(T)) \rangle,
  \qquad\text{by Lemma \ref{rightmod}}, \\ 
  & = \langle \rho, \hat{\Theta}(\pi(\omega)\sigma)(T) \rangle \\
  & =\langle \rho, \hat{\Theta}(\pi(\omega))(\hat{\Theta}(\sigma)(T)) \rangle \\
  & = \langle \rho, \hat{\Theta}(\sigma) (T) \cdot \omega \rangle, \qquad
  \text{again by Lemma \ref{rightmod}}, \\
  & = \langle \omega \bullet \rho, \hat{\Theta}(\sigma)(T) \rangle \\
  & = \langle \hat{\Theta}(\sigma)_\ast(\omega \bullet \rho), T \rangle.
  \qquad (T \in {\cal B}(L^2(G))),
  \end{split}
\]
so that
\[
  \omega \bullet \hat{\Theta}(\sigma)_\ast(\rho) = 
  \hat{\Theta}(\sigma)_\ast(\omega \bullet \rho)
\]
and therefore
\begin{equation} \label{ofi}
  \omega \bullet (\hat{\Theta}(\sigma)_\ast(\rho) - \rho) =
  \omega \bullet \hat{\Theta}(\sigma)_\ast(\rho) - \omega \bullet \rho =
  \hat{\Theta}(\sigma)_\ast(\omega \bullet \rho) - \omega \bullet \rho
  \in \tilde{I}_\sigma
\end{equation}
Since $\omega, \rho \in {\cal T}(L^2(G))$ were arbitrary, it follows that
$\tilde{I}_\sigma$ is a left ideal as claimed.
\par
To see that $\tilde{I}_\sigma$ is also a right ideal, let $\omega, \rho
\in {\cal T}(L^2(G))$, so that 
\[
  \begin{split}
  \langle (\hat{\Theta}(\sigma)_\ast(\rho)) \bullet \omega, T \rangle 
  & = \langle (\hat{\Theta}(\sigma)_\ast(\rho)), \omega \cdot T \rangle \\
  & = \langle \rho, \hat{\Theta}(\sigma)(\omega \cdot T) \rangle \\
  & = \langle \rho, \omega \cdot \hat{\Theta}(\sigma)(T) \rangle,
    \qquad\text{by Lemma \ref{leftmod}}, \\
  & = \langle \rho \bullet \omega, \hat{\Theta}(\sigma)(T) \rangle \\
  & = \langle \hat{\Theta}(\sigma)_\ast(\rho\bullet\omega), T \rangle 
  \qquad (T \in {\cal B}(L^2(G)))
  \end{split}
\]
and thus
\[
  (\hat{\Theta}(\sigma)_\ast(\rho)) \bullet \omega =  
  \hat{\Theta}(\sigma)_\ast(\rho\bullet\omega).
\]
A calculation similar to (\ref{ofi}) then lets us conclude that
$\tilde{I}_\sigma$ is a right ideal, too.
\end{proof}
\begin{remark}
The above result mirrors both \cite[Definition 3.2.1]{CL} and 
\cite[Proposition 3.3]{JN} in our setting. 
\end{remark}
\par
Finally, we consider another ideal of $({\cal T}(L^2(G)), \bullet)$.
\begin{proposition}
Let $G$ be a locally compact group, and let $\sigma \in \Mcb(A(G))$. Then
\[
  L^\infty(G)_\perp := \{ \omega \in {\cal T}(L^2(G)) : 
  \text{$\langle \omega, M_\phi \rangle = 0$ for all $\phi \in L^\infty(G)$} \} 
\]
is a two-sided ideal in $({\cal T}(L^2(G)), \bullet)$ which is contained in
the \emph{augmentation ideal}
\[
  {\cal T}_0(L^2(G)) := \{ \omega \in {\cal T}(L^2(G)) : 
                           \langle \omega, 1 \rangle = 0 \}
\]
and, if $\sigma(e) = 1$, contains $\tilde{I}_\sigma$.
\end{proposition}
\begin{proof}
Trivially, $L^\infty(G)_\perp \subset {\cal T}_0(L^2(G))$ holds.
\par
Moreover, if $\sigma(e) = 1$, then
\[
  \langle \hat{\Theta}(\sigma)_\ast(\rho) - \rho, M_\phi \rangle
  = \langle \rho,  \hat{\Theta}(\sigma)(M_\phi) \rangle - 
    \langle \rho, M_\phi \rangle
  = \langle \rho, M_\phi \rangle - \langle \rho, M_\phi \rangle
\]
holds for all $\rho \in {\cal T}(L^2(G))$ and $\phi \in L^\infty(G)$,
so that $\tilde{I}_\sigma \subset L^\infty(G)_\perp$.
\par
It remains to be shown that $L^\infty(G)_\perp$ is indeed an ideal of
${\cal T}(L^2(G))$.
\par
Let $\rho , \omega \in {\cal T}(L^2(G))$. 
Then we see that
\begin{equation} \label{yup}
  \begin{split} 
  \langle \rho \bullet \omega, M_\phi \rangle & = 
  \langle \omega , M_\phi \cdot \rho \rangle \\
  & = \langle \omega, \hat{\Theta}(\pi(\rho))(M_\phi) \rangle,
  \qquad\text{by Lemma \ref{rightmod}}, \\
  & = \langle \omega, M_\phi \hat{\Theta}(\pi(\rho))(1) \rangle \\
  & = \langle \omega, M_\phi (\pi(\rho) \cdot 1) \rangle \\
  & = \langle \rho, 1 \rangle \langle \omega, M_\phi \rangle
  \qquad (\phi \in L^\infty(G))
  \end{split}
\end{equation}
holds. From (\ref{yup}), it is immediate that $L^\infty(G)_\perp$ is indeed a 
two-sided ideal of ${\cal T}(L^2(G))$.
\end{proof}
\begin{remark}
Since $L^\infty(G)_\perp$ is a two-sided ideal of $({\cal T}(L^2(G)), 
\bullet)$, the product $\bullet$ induces a product---likewise denoted by
$\bullet$---on the quotient algebra ${\cal T}(L^2(G)) / L^\infty(G)_\perp 
\cong L^1(G)$. This product, however, is \emph{not} the usual convolution
product on $L^1(G)$: from (\ref{yup}), it is clear that
\[
  f \bullet g = \langle f, 1 \rangle g \qquad (f,g \in L^1(G)).
\]
\end{remark}
\vfill
\begin{tabbing}
{\it Second author's address\/}: \= Department of Mathematical and Statistical Sciences \kill 
{\it First author's address\/}:  \> School of Mathematics and Statistics \\
                                  \> Carleton University \\
                                  \> Ottawa, Ontario \\
                                  \> Canada K1S 5B6 \\[\medskipamount]
{\it E-mail\/}:                   \> {\tt mneufang@math.carleton.ca}\\[\medskipamount]
{\it URL\/}:                      \> {\tt http://mathstat.carleton.ca/$^\sim$mneufang/}\\[\bigskipamount]
{\it Second author's address\/}: \> Department of Mathematical and Statistical Sciences \\
                                  \> University of Alberta \\
                                  \> Edmonton, Alberta \\
                                  \> Canada T6G 2G1 \\[\medskipamount]
{\it E-mail\/}:                   \> {\tt vrunde@ualberta.ca}\\[\medskipamount]
{\it URL\/}:                      \> {\tt http://www.math.ualberta.ca/$^\sim$runde/}
\end{tabbing}
\end{document}